\documentclass[12pt]{amsart}
\usepackage{latexsym}
\usepackage{amsmath,amssymb}

 \usepackage{amscd} 
   \newtheorem{theorem}{\bf Theorem}[section]

   \newtheorem{proposition}[theorem]{\bf Proposition}
   \newtheorem{lemma}[theorem]{\bf Lemma}

\theoremstyle{remark}
   \newtheorem{example}[theorem]{Example}

\newcommand{\R}{\mathbb{R}}		
\newcommand{\C}{\mathbb{C}}		
\newcommand{\Z}{\mathbb{Z}}		
\newcommand{\ZZ}{\mathbb{Z}}		
\newcommand{\Lg}{\mathfrak{g}}		
\newcommand{\Lt}{\mathfrak{t}}		
\newcommand{\OL}{\mathcal{O}_{\Lambda}}	
\newcommand{\RR}[1]{RR(#1)}		
\newcommand{\wL}{\omega_{\Lambda}}	
\newcommand{\wo}{\omega_0}		
\newcommand{\WL}{\Omega_{\Lambda}}	
\newcommand{\WLT}{\widetilde{\Omega}_{\Lambda}}	
\newcommand{\ML}{M_{\Lambda}}		
\newcommand{\E}{\mathcal{E}}		
\newcommand{\D}{\mathcal{D}^2}		
\newcommand{\F}{\mathcal{F}}		
\newcommand{\sinhgoe}{e^{\frac{1}{2}\gamma(\E)} - e^{-\frac{1}{2}\gamma(\E)}}
\newcommand{\sinhge}{e^{\frac{1}{2}(\gamma,\E)} - e^{-\frac{1}{2}(\gamma,\E)}}
\newcommand{\sinhgex}{e^{\frac{1}{2}(\gamma,X)} - e^{-\frac{1}{2}(\gamma,X)}}

\newcommand{\ggst}{\Lg_{\text{ad}}\oplus \Lg^*_{\text{ad}}}
\newcommand{\sgn}[1]{(-1)^{#1}\,}	
\renewcommand{\t}{\mathfrak{t}}

\newcommand{\eqdef}{\,{\stackrel{{\rm def}}{=} }\,}
\newcommand{\lieg}{\Lg}
\newcommand{\liet}{\Lt} 
\DeclareMathOperator{\Td}{Td}		
\DeclareMathOperator{\ch}{ch}		
\DeclareMathOperator{\res}{res}
\DeclareMathOperator{\vol}{vol}

\newcommand{\Sl}{\mathcal{S}_{\lambda}}	
\newcommand{\Sx}{S_{\lambda}}		
\newcommand{\emm}{M}
\newcommand{\Proof}{\noindent{\bf Proof:} }

\newcommand{\calf}{{\mathcal F}}
\begin{document}

 \title[Riemann-Roch numbers of reduced spaces]
{Symplectic fibrations and 
Riemann-Roch numbers of reduced spaces}
\date{\today}

\author{Mark Hamilton}
\thanks{MH was supported in part by OGSST and University of Toronto.}
\address{University of Toronto\\Toronto, ON\\Canada}
\email{umbra@math.toronto.edu}

\author{Lisa Jeffrey}
\thanks{LJ was supported by a grant from NSERC}
\address{University of Toronto\\Toronto, ON\\Canada}
\email{jeffrey@math.toronto.edu}

\thanks{
This article makes up part of the Ph.D.\ thesis of the first author,
under the supervision of the second author.}

 \begin{abstract}
In this article we give formulas for the Riemann-Roch number of a symplectic
quotient arising as the reduced space of a coadjoint orbit
$\OL$ (for $\Lambda \in \Lg^*$ close to $0$) as an evaluation of 
cohomology classes over the reduced space
at $0$. This formula exhibits the dependence of the 
Riemann-Roch number on $\Lambda$. We also express the 
formula  as a sum over the components of the fixed point set
of the maximal torus. Our proof applies to Hamiltonian 
$G$-manifolds even if they do not have a compatible 
K\"ahler structure, using the definition of quantisation 
in terms of the Spin-$\C$ Dirac operator.

 \end{abstract}

 \maketitle

\section{Introduction}

Let $(M,\omega)$ be a compact symplectic manifold possessing a Hamiltonian
 action of a compact connected simply connected
Lie group $G$, 
with moment map $\mu: M \to {\lieg^*}$
(where $\lieg$ is the Lie algebra of $G$).  
One can form the symplectic reduction 
\[ M_0 = \mu^{-1}(0)/G,\]
 or more generally 
\[ M_{\Lambda} = \mu^{-1}(\OL)/G \] 
for $\Lambda \in \Lg,$
where $\OL \subseteq \Lg^*$ is the orbit of $\Lambda$ 
under the coadjoint action.

We assume $0$ is a regular value of $\mu$, and $0$ is a regular
value of $\mu_\Lambda: M \times \OL \to \Lg^*$
where $\mu_\Lambda(m,\xi) = \mu(m) - \xi $ 
for $m \in M$ and $\xi \in \OL$.
This is equivalent to assuming that $G$ acts with finite stabilizers
on $\mu^{-1}(0)$ (resp. $\mu_\Lambda^{-1}(0)$) \cite{GS:stp},
 so under this hypothesis
$M_0$ and  $M_\Lambda$ have at worst finite quotient
singularities. We assume that $G$ acts freely on 
$\mu_\Lambda^{-1}(0)$ and $\mu^{-1}(0)$, so that $M_\Lambda$ is a smooth 
symplectic manifold.
  Denote the symplectic form
on $M_{\Lambda}$ by $\omega_{\Lambda}$.

Let $L$ be a complex line bundle over $\ML$ with a connection whose
curvature is equal to $\wL$, called a \emph{prequantum line bundle.}
If $M$ has a complex structure compatible with the symplectic
structure (in other words $M$ is K\"ahler), then
the \emph{quantisation} of $\ML$ is defined as 
the virtual vector space
\begin{equation}
{\mathcal Q} (L^k) =  H^0(M, L^k) - H^1(M, L^k) + H^2(M, L^k) - \dots 
\end{equation}
where $H^j(M, L^k)$ is the $j$-th Dolbeault cohomology of $M$ with 
coefficients in $L^k$.
(When $k$  is very large, all the $H^j(M, L^k) = 0 $ for
$j \ne 0$, so the quantisation $\mathcal Q$ is simply a vector space.)
The dimension of the quantisation is given by the Riemann-Roch number.
The formula for the  Riemann-Roch number in terms of characteristic
classes is given below at (\ref{e:rr}).

As has been observed by Duistermaat \cite{Duist},
Guillemin \cite{G:quant}  and Vergne \cite{V:quant}, even when 
$M$ is not K\"ahler (but is equipped only with an almost complex
structure compatible with the symplectic structure -- such 
almost complex structures always exist, as explained in the above
references) one can still define the quantisation using an
elliptic complex given by the Spin-$\C$ Dirac operator.
In this more general situation,
the dimension of the quantisation is still given by the Riemann-Roch
number, the formula for which
is still given in terms of characteristic classes
by (\ref{e:rr}) below
(see \cite{Duist}, Proposition 13.1). The quantisation using the 
spin-$\C$ Dirac operator has been extensively 
studied by Meinrenken \cite{Mein1,Mein2}. 

In this article we give formulas for the Riemann-Roch number of
$L^k$ over the reduced space $\ML$
(for $\Lambda \in \Lg^*$ close to $0$) as an evaluation of 
cohomology classes over the reduced space $M_0$.
This formula exhibits the dependence of the 
Riemann-Roch number on $\Lambda$. We also express the 
formula  as a sum over the components of the fixed point set
of the maximal torus. In Section 2.1 we give a simple
proof for the K\"ahler case, which was suggested by the referee.
In Section 2.2 we treat the non-K\"ahler case.

\section{Symplectic fibrations}

It is a standard result (see for example \cite{GS:stp}) that for 
$\Lambda$ in a neighbourhood of 0 in $\Lg^*$, we have a fibration
\begin{equation} \label{e:fibn}
\begin{CD}
	\OL	@>>>	\ML\\
	@.		@VV\pi V\\
	@.		M_0
\end{CD}
\end{equation}
If $M$ is K\"ahler and 
the $G$ action preserves the 
K\"ahler structure,
  then $M_\Lambda$ and $M_0$ are also K\"ahler and (\ref{e:fibn}) is a 
fibration of K\"ahler manifolds.

Let $L$ be a line bundle over $M_{\Lambda}$ 
 with Chern character equal to $e^{\omega_{\Lambda}}$, 
(for example a prequantum line bundle) and let $k \in \Z$.
The {\em Riemann-Roch number} of $L^k$ is then
\begin{equation} \label{e:rr} \RR{\ML,L^k}=\int_{\ML} \ch(L^k)\Td(\ML) 
\end{equation}
where $\Td(\ML)$ means $\Td(T\ML)$.  
The goal of this section is to express the Riemann-Roch number
(\ref{e:rr}) 
 as far as
 possible using terms defined on $M_0$.

\newcommand{\veeb}{{\mathcal{ V} }}

Let 
\begin{equation}  \label{e:lamdef}
\lambda = k\Lambda,
\end{equation}
 where $k$ is a positive integer.
We require that $\lambda$ lie in the weight lattice
$\Lambda^W\subset \liet^*$, which is the dual of the 
integer lattice $\Lambda^I \subset \liet$ (the kernel of the 
exponential map $\liet \to T$). We do not require that $\Lambda \in 
\Lambda^W.$

\subsection{The K\"ahler case}

When $M_0$ and $M_\Lambda$ are K\"ahler, there is
the following  straightforward 
proof (which was pointed out by the referee).
Let $V(\lambda)^*$ be  the irreducible representation of $G$
with lowest weight $-\lambda$  (the dual of the irreducible representation 
with highest weight $\lambda$).
By using the principal $G$-bundle
$p_0: \mu^{-1}(0) \to M_0$, this
representation yields a vector bundle
$\veeb(\lambda)^*$ on $M_0$. 
We introduce a line bundle $L_0$ on $M_0$ for which 
$c_1 (L_0) = [\omega_0]$, where 
$\omega_0$ is the K\"ahler form of $M_0$ and 
$[\omega_0]$ is its de Rham cohomology class.

If we let $G_\Lambda$ denote the stabilizer of $\Lambda$ under the
adjoint action, then 
\begin{equation} \label{e:6}
L^k \cong (\pi^* L_0^k) \otimes {\mathcal L}_{- \lambda}
\end{equation}
where ${\mathcal L}_{- \lambda}$ is the complex line bundle associated
with the principal $G_\Lambda$-bundle
$$ p_\Lambda: \mu^{-1}(\Lambda) \to \mu^{-1}(\Lambda)/G_\Lambda
 \cong M_\Lambda$$
and with the complex representation of $G_\Lambda$ of dimension $1$ and
weight $-\lambda$.
This is true because the first Chern class of $L_k$ is equal 
to 
$L=\pi^* [\wo] + [\WLT]$, where $\WLT$ is
a form on $\ML$ which, when restricted to a fiber, is the 
Kirillov-Kostant-Souriau
 form 
$\WL$ on the coadjoint orbit $\OL$.
Furthermore, 
$$ \mu^{-1}(\Lambda)/G_\Lambda \cong \mu^{-1}(0)/G_\Lambda$$
for $\Lambda$ close to $0$.
This yields the  fibration 
$\pi$ given in (\ref{e:fibn}).
By the Borel-Weil-Bott theorem \cite{Bott}, the pushforward of 
${\mathcal L}_{-\lambda}$ under this fibration is the 
vector bundle $\veeb(\lambda)^*$, and all 
higher direct images vanish.

By the Grothendieck-Riemann-Roch theorem 
\cite{Harts}, we have
\begin{equation} \label{e:grr}
\RR{M_\Lambda, L^k} = \RR{M_0, L_0^k \otimes \veeb(\lambda)^*}
\end{equation}
using (\ref{e:6}) and the fact that higher direct images vanish.

\subsection{The non-K\"ahler case}

When $M_0$ and $M_\Lambda$ are not K\"ahler, we must rely 
on an explicit argument using the Chern character  and  the
Todd class.
By the Normal Form Theorem  (\cite{GS:stp}, 
Theorem 39.3 and Prop. 40.1) we can write the symplectic 
form on $\ML$ as $\wL=\pi^* \wo + \WLT$, where $\WLT$ is
a form on $\ML$ which, when restricted to a fiber, is the 
Kirillov-Kostant-Souriau
 form 
$\WL$ on the coadjoint orbit $\OL$.  Thus 
\[\ch(L^k) = e^{k\wL} =  e^{k\pi^* \wo} e^{k\WLT}. \]

\newcommand{\vertt}{{\mathcal T}}

Next, we can split
the tangent bundle of $\ML$ as $T\ML = \pi^* TM_0 \oplus \vertt$, 
where $\vertt$ is the vertical bundle whose fiber over a point $x\in \ML$ 
is the tangent space to the fiber of $\pi$ over $x$.
 Since the Todd class is multiplicative, 
$\Td(\ML) = \Td(T\ML) = \pi^* \Td(TM_0)\Td(\vertt).$

Combining with the expression for $\ch(L^k)$, we have 
\begin{equation}\label{rr:split}
 \RR{\ML,L^k} = \int_{\ML} \ch(L^k) \Td(\ML) 
= \int_{\ML} e^{k\pi^* \wo}\,\pi^*\!\Td(M_0) e^{k\WLT}\Td(\vertt). 
\end{equation}
The product of the 
first two factors $e^{k\pi^* \wo}\pi^*\Td(M_0)$
is written in terms of objects defined solely on the base $M_0$,
and so we turn our attention to the other factors
$e^{k\WLT}\Td(\vertt)$.
Our strategy will be to integrate 
 over the fiber of $\pi$, and be left with  an integral
 over only the base $M_0$.  
Now ${\vertt}$ is a bundle over $\ML$, and so $\Td(\vertt)
 \in H^*(\ML)$.  
If $\imath_x \colon \OL \hookrightarrow \ML$ 
is the inclusion map from 
$\OL$ to  the fiber over $x
\in M_0$, then $\imath^*_x(\vertt)\cong T\OL$.

\newcommand{\univbun}{ E}

Suppose $M$ is equipped with a 
complex vector bundle ${\mathcal{V}}$ 
(with fiber $\liet \otimes \C$) with an action of 
$G$ compatible with the action on $M$. 
Then ${\mathcal{V}}$ descends to a vector bundle $E$ on $M_0$.
The characteristic classes of $E$ come from the invariant
 polynomials on $\lieg$ via the Kirwan map.
(The Kirwan map  $\kappa: H^*_G(M) \to H^*(M_0)$ is the 
composition of the restriction map $r: H^*_G(M) \to 
H^*_G(\mu^{-1}(0))$ with the isomorphism 
$H^*_G(\mu^{-1}(0)) \cong H^*(\mu^{-1}(0)/G)$, which is 
valid when $0$ is a regular value for $\mu$.
See \cite{Ki:thesis}.)
We assume cohomology with rational, real or complex
coefficients.

We assume the 
 vector bundle $\univbun$  over $M_0$ 
has the 
property that its pullback to $M_\Lambda$ splits as the 
direct sum of a collection of line bundles $L_i$
(in other words $M_\Lambda$ is a splitting manifold for 
$\univbun$ over 
$M_0$). We then define $e_i \in H^2(M_\Lambda)$ by 
$e_i = c_1(L_i) .$ 
The characteristic class of $\univbun$ 
associated to an invariant 
polynomial $\tau \in S(\liet^*)^W$ is then given by 
$$c_\tau(\univbun) = \tau(e_1, \dots, e_\ell)$$
(where $\ell$ is the rank of $T$).
For example, if $G=U(n)$ the invariant polynomials are
generated by the elementary symmetric polynomials \cite{MS}.
The motivating
 example (the case treated in \cite{PB}) is
the case where 
$M_0$ is  the moduli space $M(n,d)$ of semistable holomorphic 
vector bundles of rank $n$ and degree $d$ over a Riemann surface
(when $n$ and $d$ are two coprime positive
integers), and
$M_\Lambda$ is the corresponding moduli space of parabolic bundles. 
In this case the vector bundle $\univbun$ is the universal bundle
(see \cite{AB}).

By Section 
14   in \cite{Hirz}, 
\[ \Td(\OL) = 
\prod_{\gamma > 0} 
	\frac{\gamma(\E)}{1-e^{-\gamma(\E)}}
=\prod_{\gamma > 0} 
	\frac{\gamma(\E) e^{\frac{1}{2} \gamma(\E)}}{(\sinhgoe)}
\]
where $\gamma$ are the roots of $G$, and $\E=(e_1,\ldots,e_\ell)$
$\in H^2(M_\Lambda) \otimes \R^\ell$.

\vspace{1ex}

For example, if $G =U(n)$,
under these hypotheses we obtain that the $j$-th Chern class is 
$$ c_j (\univbun) = \kappa (\{\tau_j\})$$
where $\tau_j$  (the $j$-th elementary symmetric polynomial)
is regarded as an element of $H^*_G({\rm pt}) = 
S(\lieg^*)^G$ and $\kappa: H^*_G(\emm) \to H^*(M_0)$ is the Kirwan map.

We introduce a basis $\hat{u_i}, i = 1, \dots, \ell$ for the 
integer lattice $\Lambda^I$ of $G$
(where $\ell$ is the rank of $G$).
This enables us to define elements
$e_j \in H^2 (M_\Lambda), j = 1, \dots, \ell $
satisfying
$c_1(L_j) = e_j$, where 
$e_j$ restricts on the fibers of $\pi$ to the generator 
$\alpha_j$ (for $ j = 1, \dots, \ell)$ of $H^2(G/T, \ZZ) \cong
H^1 (T,  \ZZ)$  corresponding to the $j$-th fundamental 
weight of $G$ (an element of ${\rm Hom}(T, U(1))$ 
 $\cong H^1(T, \ZZ))$).
Using this notation, we have

\begin{lemma}
Let $\Lambda = \sum_{i=1}^\ell \Lambda_i \hat{u_i}$.  Then 
the standard Kirillov-Kostant symplectic form $\WL$ on $\OL$ 
is given by
\[ \WL = \sum_{j=1}^\ell \Lambda_j \alpha_j,\]
where the $\hat{u_i}$ and $\alpha_j$ are as defined above.
\end{lemma}
\Proof This is a standard result (see for instance \cite{BGV}, 
Lemma 7.22).  \hfill $\Box$

\vspace{1ex}

The roots $\gamma$ also lie in  the weight lattice
$\Lambda^W$.
Writing the pairing of $\Lt^*$ and $\Lt$ as 
as $(\cdot,\cdot)$ 
we have (still when restricted to the fiber)
\[
\begin{split}
e^{k\WLT}\Td(\OL) 
& = e^{\sum \lambda_i e_i} 
	\prod_{\gamma >0} e^{\frac{1}{2}(\gamma, \E)}
\left [	\prod_{\gamma >0}\frac{(\gamma,\E)}{(\sinhge)} \right ] \\
& = e^{(\lambda,\E)} e^{\frac{1}{2}\sum\limits_{\gamma >0}(\gamma, \E)}
	\prod_{\gamma >0}\frac{(\gamma,\E)}{(\sinhge)}
\end{split}
\]

\begin{equation}\label{Wact}
 = e^{(\lambda + \rho, \E)}
	\prod_{\gamma >0}\frac{(\gamma,\E)}{(\sinhge)}
\end{equation}
\newcommand{\Wactt}{ e^{(\lambda + \rho, \E)}
	\prod_{\gamma >0}\frac{(\gamma,\E)}{(\sinhge)} }
where $\rho$ is half the sum of the positive roots.
Notice that 
\begin{equation} 
e^{k\WLT}\Td(\OL) [\OL] = \RR{\OL, L^k}
\end{equation}
which equals
$ \dim V_\lambda$ by the Bott-Borel-Weil theorem
(see \cite{Bott} or \cite{PS}),
where $V_\lambda$ is a representation of $G$ with 
highest weight $\lambda$ (we have assumed that 
$\lambda$ is in the fundamental Weyl chamber).

After evaluating on the fundamental cycle of 
 $M_\Lambda$, the equation (\ref{Wact}) equals 
\begin{equation} \label{e:top4}
 \frac{1}{\lvert W \rvert} \left( \, \sum_{\sigma \in W} \sgn{\sigma}
	\frac{e^{(\sigma(\lambda + \rho) , \E)}}{\prod_{\gamma>0}(\sinhge)}
  \right) \prod_{\gamma>0} (\gamma,\E) [M_\Lambda].
\end{equation}
The expression
in brackets in \eqref{e:top4}
 is unchanged under the action of the 
Weyl group  $\E \mapsto w\E$.

Let $\{\tau_r\}$ ($r = 1, \dots, n(G)$) be a set of generators for the 
ring $S(\liet^*)^W$ of 
Weyl invariant polynomials on $\liet$, where $n(G)$ is the number 
of generators.
By Proposition 3.6 in \cite{PB}, $\tau_r(e_1, \dots, e_\ell)) = \pi^*a_r$
for $a_r$ a class on $M_0$.
Therefore the factor (\ref{e:top4})
 can be written as a function of 
the invariant polynomials $\tau_r$ applied to the $e_j$'s,
\begin{equation} \label{e:Slambda}
 \Sl \bigl(  \tau_r(e_1, \dots, e_\ell)  \bigr) 
\eqdef \sum_{\sigma \in W} \sgn{\sigma}
	\frac{e^{(\sigma(\lambda + \rho) , \E)}}{\prod_{\gamma>0}(\sinhge)}
\end{equation}
We see that $\Sl$ is actually a polynomial
in the $e_j$. 
Because $e_j^{N+1} = 0 $ where
$2N = \dim_{\R} M_\Lambda$,
it follows immediately  that after evaluating
on the fundamental class of $M_0$,
\begin{theorem} \label{t:degreef} 
 $\Sl$
is a polynomial in $\lambda$ of degree $\le N$.
\end{theorem}

We can replace the term $e^{k\WLT}\Td(V)$ in the integral 
with $$\frac{1}{|W|}
\Sl(\pi^*a_1, \dots, \pi^* a_{n(G)}) \prod_{\gamma>0}(\gamma,\E)$$
to get
\begin{equation} \label{e:rrk}
\RR{\ML,L^k} = \frac{1}{|W|} \int_{\ML} e^{k\pi^*\wo} \pi^*\Td(M_0) 
	 \Sl (\pi^* a_1, \dots, \pi^* a_{n(G)})
	\prod_{\gamma>0} (\gamma,\E).
\end{equation}
All of the factors in the integral except for $\prod (\gamma,\E)$
are constant on each fiber, and so we have
\begin{equation} \label{e:rr2}
\RR{\ML,L^k} = \frac{1}{|W|}\int_{M_0} e^{k\wo} \Td(M_0) \Sl(a_1, 
\dots, a_{n(G)})
	\int_{\OL}\prod_{\gamma>0} (\gamma,\E).
\end{equation}
Since \[ \int_{\OL}\prod_{\gamma>0} (\gamma,\E) = |W|,\] 
we have finally
\begin{theorem} \label{t:rrfin}
The Riemann-Roch number of a symplectic fibration $M_\Lambda$ is 
given by 
\begin{equation} \label{e:rrfin}
\RR{\ML,L^k} = \int_{M_0} e^{k\wo} \Td(M_0) \Sl(a_1,\dots, a_{n(G)}).
\end{equation}
where $\Sl $ is defined at (\ref{e:Slambda}).
\end{theorem}

\begin{theorem} \label{t:asymptotics}
When $\Lambda \in \Lambda^W$ is a weight in the fundamental
Weyl chamber,
 then the limit as 
$k \to \infty $ of 
$$\frac{\RR{M_\Lambda, L^k}}{k^{N}}$$  
(where as above $2N = \dim_\R (M_\Lambda))$ is
$   \vol M_0 \dim V_{\Lambda - \rho} $. 
\end{theorem}
\Proof 
This limit is given by $\vol M_\Lambda$. The symplectic 
volume is given by
$\vol M_0 \vol \OL$.
The symplectic volume of $\OL$  is given  at \cite{BGV} (Proposition 7.26)
as
\begin{equation} \label{e:symvol}
\vol (\OL) = \frac{\prod_{\alpha > 0 } <\alpha, \Lambda> }
{\prod_{\alpha > 0 } <\alpha, \rho> }.
\end{equation}
This gives the value of 
$\dim V_{\Lambda-\rho}$, 
using the Weyl dimension formula \cite{DK}.
\hfill $\Box$

\begin{proposition} \label{p:2.5}
Let $\lambda \in \Lambda^W$ be a   weight in the 
fundamental Weyl chamber, and define $\Lambda (k) = \lambda/k$.
Let $N_0 = \frac{1}{2} {\rm dim} M_0$.
Then $$\lim_{k \to \infty} \frac{1}{k^{N_0} }RR(M_\Lambda, L^k) 
= ({\rm vol} M_0) (\dim V_{\lambda - \rho} ).  $$
\end{proposition}
\Proof For $X \in {\bf t }$
we introduce 
\begin{equation}
S_\lambda\bigl(  X  \bigr) 
\eqdef \sum_{\sigma \in W} \sgn{\sigma}
	\frac{e^{(\sigma(\lambda + \rho) , X)}}{\prod_{\gamma>0}(\sinhgex)}
\end{equation}
Notice that we are fixing $\lambda$ and allowing $\Lambda = \lambda/k$
to vary as $k$ varies.
Notice that by the Weyl character formula (\cite{PS}, Proposition 
14.2.2)  we have 
$$S_\lambda (X) = \chi_\lambda (\exp X) $$
where $\chi_\lambda $ is the character of the representation 
with lowest weight  $-\lambda$.
Theorem \ref{t:rrfin} gives 
the result, noting that in the limit $k \to \infty$ 
the leading order term in 
$k$ comes by integrating $(k \omega_0)^{N_0}$ 
so the factor $\Sl (a_1, \dots, a_{n(G)}) $ contributes only 
its value when all the arguments $a_i$ are 
replaced by $0$, in other words when 
the argument of $S_\lambda (X)$ is replaced by 
$X = 0 $. This is 
the value $\chi_\lambda (0)$, in other
words the dimension of $V_{\lambda - \rho}$. 
\hfill $\square$

\begin{example} When $G = SU(2)$, 
the value of 
$\Sl$ is  (recalling that $\lambda$ is a positive integer)
$$ \Sl(a_2) = \kappa(S_\lambda(X) ) $$
where 
\begin{equation} \label{e:su2}
S_\lambda(X) = \frac{1}{2} \frac{e^{(\lambda + 1) X} - e^{- (\lambda + 1) X} }
{e^X - e^{-X} } 
\end{equation} 
\begin{equation}= \frac{1}{2} (\cosh(X) + \dots + \cosh (\lambda - 1) X 
+ \cosh \lambda X ). \end{equation} 
Here we have introduced a formal variable $X$ for which
$\kappa (X^2) = a_2 \in H^4(M_0). $ 
It follows that 
$\Sl$ is a polynomial in $\lambda$ of order 
$N + 1$, because the terms which contribute from the 
Taylor expansion of  order $X^N $ are 
$ \sum_{j = 1}^\lambda j^N$ which is of order 
$\lambda^{N+1}$.
Because we are integrating over $\ML$, the
 ring of polynomials in the variable 
$X$ gets truncated by imposing the 
relation $X^{N+1} = 0 $.
In this example $S_\lambda (0) = \lambda/2. $ 
\end{example}

\section{The Jeffrey-Kirwan residue formula}

In the final two sections of this paper, we express the Riemann-Roch 
number on a symplectic fibration in terms of data at the fixed point
set of the maximal torus $T$ of $G$ (assuming $\Lambda$ is generic so
its stabilizer under the coadjoint action is $T$).
In this way, we obtain a second proof of Theorem \ref{t:rrfin}.

The residue formula \cite{nonloc}
expresses cohomology pairings on reduced spaces $M_0$ in
terms of a multi-dimensional residue of certain rational holomorphic
functions on $\mathfrak{t}$ \cite{nonloc}. 
This formula is valid provided $0$ is a regular value of the moment map.
The cohomology classes $\beta_0$ on $M_0$ are assumed to come from 
equivariant cohomology classes $\beta$ on $M$ via the Kirwan map
\cite{Ki:thesis}.
The fixed point data are
\begin{itemize}
\item the value of the moment map at a component $F$ of the fixed point set of 
the maximal torus $T$
\item the restriction of $\beta$ to the $F$
\item the equivariant Euler class $e_F$ of the normal bundle to $F$
(which involves the weights of the action of $T$ on the 
normal bundle, as well as the ordinary Chern roots of the 
normal bundle).
\end{itemize}
The residue formula takes the form
\begin{equation} \label{e:resgen}
\int_{M_0} e^{\omega_0} \beta_0 = C {\rm Res} \left (\sum_F \int_F 
\frac{e^{\omega + \mu(F)(X)} 
\beta(X)}{e_F(X)} \right ) \end{equation} 
where $C$ is a nonzero
constant, $X \in \liet\otimes \C$
 is a formal variable (in the Cartan model for 
equivariant cohomology), and ${\rm Res}$ is defined at (\ref{e:resgen}) 
below. 

We can readily identify the equivariant cohomology class giving rise to 
the Riemann-Roch number of a prequantum line bundle. The class
$e^{\omega_0}$ is the Chern character of the line bundle
over $M_0$, while
$e^{\omega + \mu(F)(\cdot)}$ is the equivariant Chern character of the 
line bundle over $M$.
The relevant class $\beta_0$ is the Todd class of $M_0$,
which arises in the image of the 
Kirwan map using the equivariant Todd class
${\rm Td}_G$ 
(see \cite{locq} and Proposition \ref{p:todd} below). 
The residue formula has been applied to studying Riemann-Roch
numbers 
in \cite{LRR} and \cite{locq}.
Here we study  the residue formula  for the Riemann-Roch number
on a symplectic fibration.

The residue formula 
applies to compact $M$ reduced by any compact group $G$.
The computation of terms in the residue formula depends on 
the choice of
a cone
$\Gamma$ in $\mathfrak{t}$, even though the result of the formula is
independent of this choice. Let $\gamma_1,\dots,\gamma_k$ be the set
of all weights that occur by the $T$ action at any of the fixed point
components. Choose some $\xi\in {\bf t}$ such that $\gamma_i(\xi)\neq 0$ for
all $i$. Let $\beta_i = \gamma_i$ if $\gamma_i(\xi)>0$ and
$\beta_i=-\gamma_i$ if $\gamma_i(\xi)<0$. Thus $\beta_i(\xi)>0$ for
all $i$. The cone $\Gamma$ is the set of all vectors in $\t$ which
behave like $\xi$:
$$
\Gamma = \{X\in\t:\ \beta_i(X)>0,\ \mbox{for all}\ i\}.
$$

\begin{theorem}[Jeffrey-Kirwan]  \label{t:residue}
Let $(M,\omega)$ be a compact symplectic
manifold with a Hamiltonian $T$ action and moment map $\Phi$, where
$T$ is a compact torus. Denote by ${\mathcal F}$ the connected 
components of the fixed point set of $T$ on $M$. Let $p$ be a regular value of $\Phi$  and
$\omega_p$ the Marsden-Weinstein reduced symplectic form on
$M_p$. Then for $\beta\in H^*_T(M)$ and $\kappa_p:H_T^*(M)\rightarrow
H^*(M_p)$ we have 
\begin{equation} \label{e:star}
\int_{M_p}\kappa_p(\beta)e^{\omega_p} = C\cdot \res^\Gamma \left(\sum_{F\in
{\mathcal F}} e^{i(\Phi(F)-p)(X)}\int_F\frac{\iota^*_F (\beta(X)
e^{\omega})}{e_F(X)}[dX]\right) 
\end{equation}
where $C$ is a non-zero constant, $X$ is a variable in
$\mathfrak{t}\otimes \C$, and $e_F(X)$ is the equivariant Euler class
of the normal bundle to $F$ in $M$. The multi-dimensional residue
$\res^\Gamma$ is defined below at (\ref{e:resdef}).
\end{theorem}

The residue can be  defined as follows (see \cite{locq} Proposition 3.4).
For $f$ a meromorphic function of one complex variable
$z$  which is of the 
form $f(z) = g(z) e^{i \lambda z} $ where $g$ is a rational function,
we define
$$\res^+_z f(z)dz =  
\sum_{b\in\C} \res(g(z)e^{i\lambda z};z=b).$$
We extend this definition  by 
linearity to linear combinations of functions
of this form.

Viewing $f$ as a meromorphic function on the Riemann sphere and observing
that the sum of all the residues  of a meromorphic 1-form
on the Riemann sphere is $0$, we observe that 
$$ \res^+_z \left ( f(z)dz\right )  = 
- {\rm Res}_{z = \infty} \left ( f(z)dz\right ). $$

If $X  \in \liet$,
define
$$h(X)=\frac{q(X)e^{i\lambda(X)}}{\prod_{j=1}^k \beta_j(X)}$$
for some polynomial function $q(X)$ of $X$
 and some $\lambda,\beta_1,\ldots,\beta_k \in {\liet}^*$.
Suppose that $\lambda$ is
 not in any proper subspace of ${\liet}^*$ spanned by a subset of
$\{\beta_1,\ldots,\beta_k\}$. Let $\Gamma$
 be any nonempty open cone in $\liet$ contained in
some connected component of 
$$\{X\in\liet : \beta_j(X) \neq 0, 1\leq j\leq k\}.$$
Then for a generic choice of coordinate system $X=(X_1,\ldots,X_l)$ 
on $\liet$ for
which  $(0,\ldots,0,1)\in \Gamma$ we have
\newcommand{\Jac}{{\rm Jac}}
\begin{equation} \label{e:resdef}
\res^{\Gamma}(h(X)[dX])=  \Jac 
\res^+_{X_1} \circ  \ldots \circ \res^+_{X_l}
\left (  h(X)dX_1\ldots dX_l\right )\end{equation}
where the variables $X_1,\ldots,X_{m-1}$ are held constant while calculating
$\res^+_{X_m}$, and $\Jac$ is the determinant of any $l\times l$ matrix whose
columns are the coordinates of an orthonormal basis of $\liet$ defining the
same orientation as the chosen coordinate system. 
We assume that if $(X_1, \dots, X_l)$ is a coordinate
system for $X \in \liet$, then
$(0, 0, \dots, 1) \in \Gamma$. We  also require an additional
technical  condition on the coordinate systems, which
is valid for almost any choice of coordinate system (see
Remark 3.5 (1) from \cite{locq}).

\section{Fixed point formulas}	
We wish to use the residue formula to calculate the integral 
in \eqref{e:rr2}.  
In (\ref{e:star})
we are interested in 
$ \beta(X) = \Td_G(M\times \OL) \Td_G^{-1}(\ggst)$
which maps to $\Td (M_\Lambda)$ under the Kirwan map.


The components of the 
fixed point set for the action of $T$ on $M \times \OL$
are of the form $ F\times \{\sigma\Lambda \} $ where $F $ is a component of the 
$T$ fixed point set of $M$ and $\sigma \in W$.
The equivariant Euler class at this fixed point is 
$e_F(X) (-1)^\sigma\prod_{\gamma> 0} \gamma(X).$ 
Thus the residue formula becomes 

$$
\int_{M_\Lambda}  e^{k \omega_\Lambda} \Td (M_\Lambda)
 = 
C  {\rm Res}\Biggl ( \sum_{\sigma \in W} \sum_{F \in \calf} (-1)^\sigma
e^{\sigma \lambda( X)}e^{k\mu_F(X)} \times
$$
\begin{equation} \label{e:restwo}
\times \frac{\D(X)  }{ \prod_{\gamma>0} \gamma(X) \Td_G(\lieg_{\rm ad} \oplus
 \lieg_{\rm ad}^*) }
 \int_{F\times \{\sigma \Lambda\}} \frac{e^{k\omega}  } 
{e_F(X) } {\Td_G(M) \Td_G (\OL)} 
  \Biggr ). \end{equation}
Here,
 $\omega$ is the symplectic form on $F$,
 $C$ is an overall constant
given in the
statement of 
Theorem \ref{t:residue},
 $\D(X)$ is the product of all the roots of $\lieg$, and
$e_F(X) \in H_T^*(M)$ is the equivariant Euler class of the normal bundle
to $F$. The  $F$ are 
components of the fixed point set $\F$ for the
action of $T$ on $M$. We use Proposition \ref{p:todd} below for the Todd class.

\begin{proposition} \label{p:todd} (\cite{LRR}, Proposition 2.1)
The formal equivariant cohomology class 
\begin{equation}\label{kappa}
 \Td_G(M) \Td_G^{-1}(\ggst) 
\end{equation}
maps to $\Td(M_0)  $
 under $\kappa$,
where $\Lg_{\text{ad}}$ (resp.  $\Lg^*_{\text{ad}}$) denotes the product
bundle $M \times \Lg$ (resp. $M \times \Lg^*$)
 with $G$ acting on $\Lg$ by the adjoint action
(resp. the  coadjoint action).
\end{proposition}
\Proof
We observe that  $\kappa \bigl( \Td_G(M) \Td_G^{-1}(\ggst) \bigr) = \Td(M_0)$. 
\hfill $\Box$

\begin{proposition} \label{p:sx}
The formal equivariant cohomology class $\Sx(X)$
maps to $\Sl (a_1, \dots, a_{n(G)})$
under $\kappa$,
where
\begin{equation} \label{e:Sdef}
	\Sx(X)=	\frac{\sum_{\sigma \in W} \sgn{\sigma}
			e^{(\sigma(\lambda + \rho) , X)}}
		{\prod_{\gamma>0}(\sinhgex)}.
\end{equation}
\end{proposition}

\Proof
The expression defining $\Sx$ is the same 
as the expression (\ref{e:Slambda}) defining $\Sl$, with the 
expression $\E$ replaced by the variable $X \in \liet$.
Thus $\Sx \in S(\liet^*)^W$ can be viewed as an equivariant 
cohomology class on $M$.  

Since $\Sx$ is symmetric under the action of the 
Weyl group, it is a function of the invariant 
polynomials $\tau_r(X)$.  Since $\kappa(\tau_r)=a_r$ by
definition of the $a_r$, 
when $\kappa$ is applied to $\Sx$, the result will be $\Sl$.  
Since $\kappa$ is a ring homomorphism, the result follows.
\hfill $\Box$

The expression (\ref{e:restwo}) is equal to 
the expression (\ref{e:top4}) which we derived in Section 2 for 
\begin{equation} \label{e:res}
	\int_{M_0} e^{k\wo} \Td(M_0) \Sl(a_1, \dots, a_{n(G)}).
\end{equation}
To see this,
we only need to evaluate (\ref{e:top4}) using
the residue formula, and use the fact that 
$$\frac{1}{\prod_{\gamma>0} \gamma(X)} 
\Td_G(\OL)|_{F \times \{\sigma \Lambda\} }
 = (-1)^\sigma \prod_{\gamma>0}\frac{1}{1 - e^{-\sigma \gamma(X) }} $$
(which appears in (\ref{e:restwo}))
is the same as 
$$\frac{e^{\sigma \rho}}{\prod_{\gamma>0} (e^{\gamma/2} - e^{-\gamma/2})} $$
(which appears  when we evaluate (\ref{e:top4})
using the residue formula).
Hence \eqref{e:res} becomes
\begin{equation}\label{RRres}
	C\res \Bigl ( \D(X) 	\Sx(X) 
\sum_{F\in \F} e^{i\mu_T(F)(X)}
	\int_F \frac{\imath_F^* \bigl( 
		\Td_G(M) \Td_G^{-1}(\ggst) e^{k\omega} \bigr)}
		{e_F(X)}
	[dX] \Bigr),
\end{equation}
which is the same as (\ref{e:restwo}).
Thus we have obtained a second proof of Theorem \ref{t:rrfin}.

\end{document}